\documentclass{amsart}

\usepackage{amssymb}
\usepackage{amsmath}
\usepackage{mathrsfs}
\usepackage[english,francais]{babel}
\usepackage[latin1]{inputenc}
\usepackage[usenames]{color}

\newcommand{\R}{\mathbb R}
\newcommand{\Hh}{\mathcal{H}}

\newcommand{\ve}{\varepsilon}

\newtheorem{theorem}{Theorem}[section]

\newtheorem{e-proposition}[theorem]{Proposition}

\newtheorem{e-definition}[theorem]{Definition\rm}

\def\og{\leavevmode\raise.3ex\hbox{$\scriptscriptstyle\langle\!\langle$~}}
\def\fg{\leavevmode\raise.3ex\hbox{~$\!\scriptscriptstyle\,\rangle\!\rangle$}}

\begin{document}
\centerline{}


\title{A Modica-Mortola approximation for the Steiner Problem}
\thanks{This work has been partially supported by the Agence Nationale de la Recherche, through the project ANR-12-BS01-0014-01 GEOMETRYA, and by The Gaspard Monge Program for Optimization and operations research (PGMO) via the project MACRO}

\author[A. Lemenant]{Antoine Lemenant}
\author[F. Santambrogio]{Filippo Santambrogio}

\address[A. Lemenant]{Universit\'e Paris-Diderot, Laboratoire Jacques-Louis Lions}
\email{lemenant@ljll.univ-paris-diderot.fr}

\address[F. Santambrogio]{Universit\'e Paris-Sud, Laboratoire de Math\'ematiques d'Orsay}
\email{santambrogio@math.u-psud.fr}

\maketitle
\selectlanguage{english}




\begin{abstract}
\selectlanguage{english}
In this note we present a way to approximate the Steiner problem by a family of elliptic  energies
of Modica-Mortola type, with  an additional  term relying on the weighted geodesic distance which takes care of the connexity constraint. 
\vskip 0.5\baselineskip

\selectlanguage{francais}
\noindent{\bf R\'esum\'e} \vskip 0.5\baselineskip \noindent

Dans cette note nous pr\'esentons une m\'ethode d'approximation du probl\`eme de Steiner par une famille de fonctionnelles de type Modica-Mortola, avec un terme additionnel bas\'e sur une distance g\'eod\'esique \`a poids, pour prendre en compte la contrainte de connexit\'e. 
\medskip

\noindent {\bf Titre: Une approximation à la Modica-Mortola pour le probl\`eme de Steiner.}

 \end{abstract}

\selectlanguage{francais}
\section*{Version fran\c{c}aise abr\'eg\'ee}

Le  probl\`eme bien connu dit ``de Steiner" consiste \`a trouver un compact connexe de longueur minimal qui contient certains points du plan donn\'es au d\'epart, en nombre fini. L'ensemble minimal est alors un arbre fini constitu\'e de segments qui peuvent se joindre par nombre de 3 uniquement, formant des angles de 120° \cite{gilbpoll,ps2009}. L'un des aspects qui a rendu ce probl\`eme si c\'el\`ebre réside dans sa complexité de calcul, malgré une formulation simple en apparence, faisant partie de la liste des 21 probl\`emes NP-complets de Karp \cite{Karp} (le temps polynomial \'etant \'evalu\'e par rapport au nombre de points).

Dans cette note nous proposons une m\'ethode susceptible de donner lieu \`a des solutions approch\'ees du Probl\`eme de Steiner. La strat\'egie repose sur  l'emploi de fonctionnelles de type elliptique \`a la mani\`ere de Modica-Mortola \cite{MM}, comme l'ont fait d'autres auteurs auparavant concernant des probl\`emes li\'es au p\'erim\`etre ou longueur d'un ferm\'e \cite{at,OudetKelvin,sant,OudetSantam,alr,dmi}. La nouveaut\'e dans notre approche est l'ajout d'un terme permettant de g\'erer la contrainte de connexit\'e sur l'ensemble \`a minimiser.

Ce nouveau terme fait intervenir la fonction distance pond\'er\'ee $d_\varphi$, d\'efinie en \eqref{defDphi}. Cette fonction peut être calcul\'ee num\'eriquement sur une grille par une m\'ethode, dite {\it fast-marching} \cite{sethian-fastmarching}, qui a \'et\'e r\'ecemment am\'elior\'ee dans \cite{bcps} permettant le calcul \`a la fois de $d_\varphi$ et de son gradient par rapport \`a $\varphi$. La fonctionnelle approximante que nous proposons est la suivante
\begin{equation*}
S_\varepsilon(\varphi):=\frac{1}{4\varepsilon}\int_{\Omega}(1-\varphi)^2dx + \varepsilon\int_{\Omega}\|\nabla \varphi\|^2+\frac{1}{\sqrt{\varepsilon}}\sum_{i = 1}^{N}d_\varphi(x_i,x_1). \label{defSeps0}
\end{equation*}
Notre r\'esultat principal stipule qu'étant donné une suite de minimiseurs $\varphi_\ve$ de $S_\varepsilon$ et la suite de fonctions $d_{\varphi_\ve}(\cdot,x_1)$ associée, ces fonctions convergent \`a une sous-suite pr\`es vers une fonction $d$, dont l'ensemble de niveau $\{d=0\}$ est un minimiseur du Probl\`eme de Steiner associ\'es aux points $\{x_i\}$ (Th\'eor\`eme \ref{SteinerTheorem}).

Les deux premiers termes de la fonctionnelle rappellent la fonctionnelle de Modica-Mortola. Le point essentiellement nouveau demeure dans l'ajout du troisi\`eme terme de la fonctionnelle, bas\'e sur le fait suivant: si $\sum_{i = 1}^{N}d_\varphi(x_i,x_1)=0$, alors l'ensemble $\{d_\varphi =0\}$ doit être connexe par arcs et contenir les $\{x_i\}$.

Dans l'article \cite{bls}, \'ecrit conjointement avec M. Bonnivard, nous utilisons cette technique pour approcher \'egalement certaines variantes du probl\`eme de Steiner, comme par exemple la fonctionnelle de distance moyenne. On peut y trouver des preuves d\'etaill\'ees ainsi que des simulations num\'eriques.

Les m\'ethodes num\'eriques envisag\'ees, inspir\'ees par le travail d'\'E. Oudet dans \cite{OudetKelvin,OudetSantam}, se basent sur une m\'ethode de gradient appliqu\'ee \`a chaque fonctionnelle $S_\ve$ (qui est convexe pour $\ve$ grand), en diminuant par \'etapes la valeur de $\ve$ et prenant comme initialisation \`a chaque \'etape le point de minimum approch\'e trouv\'e \`a l'\'etape pr\'ec\'edente. Cela ne garantie pas de converger vers un minimum global, mais permet en g\'en\'eral de choisir un ``bon'' minimum local. 

La preuve du th\'eor\`eme d\'ecrit plus haut est de type $\Gamma$-convergence. Plus pr\'ecis\'e-ment, diff\'erem\-ment de ce qui a \'et\'e fait dans  \cite{OudetKelvin,OudetSantam} ainsi que dans les autres cas \'etudi\'es dans  \cite{bls}, il n'est pas possible de mani\`ere \'evidente d'exprimer notre r\'esultat sous la forme d'un \'enonc\'e de $\Gamma-$convergence d'une suite de fonctionnelles vers une autre. Cependant, la d\'emonstration en suit le m\^eme sch\'ema. La $\Gamma$-limsup d\'ecoule de techniques classiques que l'on peut trouver dans \cite{at}. En revanche la $\Gamma$-liminf est plus d\'elicate. L'un des points difficiles \`a montrer est la rectifiabilit\'e d'une limite Hausdorff d'ensembles de niveaux de fonctions $d_{\varphi_\varepsilon}$ associ\'ees \`a des $\varphi_{\varepsilon}$ d'\'energies uniform\'ement born\'ees. L'argument original de Modica-Mortola \cite{MM} est bien s\^ur essentiel, mais de nouvelles techniques n\'ecessitent d'être introduites.




\selectlanguage{english}
\section{Introduction}
\label{}

Given a finite number of points $D:=\{x_i\}_{i=1,\dots,N}\subset \Omega \subset  \R^2$, the so-called Steiner problem consists in solving
\begin{eqnarray}
\min\big \{\Hh^1(K) \quad ; \; K\subset \R^2 \text{ compact, connected, and containing } D \big\}. \label{Steiner00}
\end{eqnarray}

Here, $\Hh^1(K)$ stands for the one-dimensional Hausdorff measure of $K$. It is  known that minimizers for \eqref{Steiner00} do exist, need not to be unique, and are  trees composed by a finite number of segments joining with only triple junctions at 120°, whereas  computing a minimizer is very hard (some versions of the Steiner Problem belong to the original list of NP-complete problems by Karp, \cite{Karp}). We refer for instance to \cite{gilbpoll} for a history of the problem and to \cite{ps2009} for recent mathematical results about it.



In this note we propose a way to approximate the problem, and we prove convergence to an exact solution as some parameter $\varepsilon$ goes to zero. Our strategy is to approximate the length by an elliptic energy of Modica-Mortola \cite{MM} type. This strategy was pursued before by many authors for similar problems involving the perimeter  or the length of a closed set (see e.g. \cite{at,OudetKelvin,sant,OudetSantam,alr,dmi}), but the novelty here is that we are able to add a term taking care of the connexity constraint. This term relies on the weighted geodesic distance  $d_\varphi$, defined as follows. Given $\Omega \subset \R^2$, for any non-negative function $\varphi \in C^0(\overline{\Omega})$, we define the corresponding weighted geodesic distance through

\begin{equation}
d_\varphi(x,y):= \inf \left\{ \int_{\gamma} \varphi(x) d\mathcal{H}^1(x) ; \; \gamma \text{ curve in $\Omega$ connecting } x \text{ and } y \right\}. \label{defDphi}
\end{equation}

Given a function $\varphi$ and a point $x_1$, the distance $d_\varphi(\cdot,x_1)$ can be treated numerically by the so-called {\it fast-marching} method  \cite{sethian-fastmarching} since it is a solution of $\|\nabla u \|=\varphi$ with $u(x_1)=0$ in the viscosity sense. A recent improvement of this algorithm (see \cite{bcps}) is now able to compute at the same time $d_\varphi$ and its gradient with respect to $\varphi$, which is useful every time one needs to optimize w.r.t. $\varphi$ a functional involving $d_\varphi$. Our proposal to approximate the problem \eqref{Steiner00} is then to minimize
$$S_\varepsilon(\varphi):=\frac{1}{4\varepsilon}\int_{\Omega}(1-\varphi)^2dx + \varepsilon\int_{\Omega}\|\nabla \varphi\|^2+\frac{1}{\sqrt{\varepsilon}}\sum_{i = 1}^{N}d_\varphi(x_i,x_1),$$
among all functions $\varphi \in \mathcal{A}:=H^{1}(\Omega)\cap C^0(\overline{\Omega}) \cap\{\ve\leq \varphi\leq 1 \text{ and }\varphi=1 \text{ on } \partial \Omega\}$.

The first two terms are a simple variant of the standard Modica-Mortola functional, already used in  \cite{at}: as $\ve\to 0$, they force $\varphi$ to tend to $1$ a.e. and pay the transition between the value $1$ and the value $\ve$ by means of the length of the transition set, while the last term tends to enforce connectedness. The key point is that whenever $\sum_{i = 1}^{N}d_\varphi(x_i,x_1)=0$, the set $\{d_\varphi =0\}$ must be path-connected, must contain all the points $\{x_i\}$, and the path connecting them inside this set are such that $\varphi=0$ $\mathcal{H}^1-$a.e. on them.


In the paper \cite{bls}, the authors together with M. Bonnivard used this idea to approximate some variant of  the Steiner Problem,  as the Average distance and $p$-Compliance problem.  One can find therein detailed proofs and numerical experiments.
The main idea for numerics is based on the work by \'E. Oudet in \cite{OudetKelvin,OudetSantam}: for every $\ve$ one can run a gradient descent for $S_\ve$ (which is convex for large $\ve$), and a candidate minimizer for the limit problem is obtained by  reducing at each step the value $\ve$ and initializing the gradient with the critical point obtained at the previous step. There is no guarantee that this converges to a global minimum, but at least a ``well-chosen'' local minimum is chosen.

{\bf Existence of minimizers for $S_\ve$.} The existence of minimizers for the functional $S_\ve$ is a delicate matter. This depends on the fact that $H^1$ does not inject into $C^0$ and on the behavior of the map $\varphi\mapsto d_\varphi$. First, notice that we only restricted our attention to $\varphi\in C^0(\overline{\Omega})$ for the sake of simplicity. Indeed, it is possible to define $d_\varphi$ as a continuous function as soon as $\varphi\in L^p$ for an exponent $p$ larger than the dimension (here, $p>2$, see \cite{CarJimSan}). The difficult question is which kind of convergence on $\varphi$ provides pointwise convergence for $d_\varphi$. If one wanted upper semi-continuity of the map $\varphi\mapsto d_\varphi(x,x_1)$ (for fixed $x$ and $x_1$), this would be easy, thanks to the concave behavior of $d_\varphi$, and any kind of weak convergence would be enough. Yet, in this case we would like lower semi-continuity, which is more delicate. An easy result is the following: if $\varphi_n\to\varphi$ uniformly and a uniform lower bound $\varphi_n\geq c>0$ holds, then $d_{\varphi_n}(x,x_1)\to d_{\varphi}(x,x_1)$. Counterexamples are known if the lower bound is omitted. On the contrary, replacing the uniform convergence with a weak $H^1$ convergence (which would be natural in the minimization of $S_\ve$) is a delicate matter (by the way, the continuity seems to be true and it is not known whether the lower bound is necessary or not).

For the sake of our paper, one could consider adding an extra term of the form $\ve^{10}\int ||\nabla\varphi||^p$ with $p>2$, which enforces continuity and uniform convergence, or just think that the results are given ``provided a minimizer exist''. From the point of view of the approximation result and of the numerical applications this is not crucial.

\section{The main result }

\begin{theorem} \label{SteinerTheorem} Let $\Omega$ be a bounded open convex set containing the convex hull of the $\{x_i\}$. For all $\varepsilon>0$ let  $\varphi_{\varepsilon}$ be a minimizer of $S_{\varepsilon}$ among all $\varphi \in \mathcal{A}$. Consider the sequence of functions $d_{\varphi_\ve}(\,\cdot\,,x_1)$, which are $1-Lipshitz$ and converge, up to subsequences, to a certain function $d$. Then the set $K:=\{d=0\}$ is compact, connected and is a solution to the Steiner Problem \eqref{Steiner00}.
\end{theorem}

{\bf Proof.} We first extract a subsequence such that the sequence of $1$-Lipschitz functions $d_{\varphi_{\varepsilon_n}}(x,x_1)$ converges uniformly to some function $d(x)$. It is easy to see that the set  $K=\{d(x)=0\}$,  is a compact and connected set as a Hausdorff limit of sub level sets of $d_{\varphi_{\varepsilon_n}}(\;\cdot\;,x_1)$, which are all compact connected sets.

Let now $K'$ be any competitor in the Steiner Problem, that we can assume contained in $\Omega$ . By using a variant of  \cite[Theorem 3.1.]{at}, it is not difficult to construct a sequence of functions $\psi_\varepsilon \in H^1(\Omega)\cap C^0(\overline{\Omega})$, satisfying $\ve\leq \psi_\ve\leq 1$, $\psi_\varepsilon=1$ on $\partial \Omega$ and $\limsup_{n}S_{\varepsilon_n}(\psi_{\varepsilon_n})\leq \Hh^1(K')$.
In particular, following the construction of \cite{at} it is easy to make the last term $\frac{1}{\sqrt{\varepsilon}}\sum_{i = 1}^{N}d_{\varphi_\varepsilon}(x_i,x_1)$ tend to zero since $\varphi_\varepsilon$ is  very small close to $K'$ (a careful look at the proof reveals that $1/\sqrt{\varepsilon}$ is needed in front of this term, or any other coefficient of the form $o((\ve\ln\ve)^{-1})$).

On the other hand it is clear from the minimizing property of $\varphi_\varepsilon$ that 
\begin{equation}
\liminf_nS_{\varepsilon_n}(\varphi_{\varepsilon_n})\leq \limsup_{n}S_{\varepsilon_n}(\psi_{\varepsilon_n}),
\end{equation}
 thus the proof will be finished provided that we show the following claim
\begin{eqnarray}
\Hh^1(K)\leq \liminf_{n\to +\infty} S_{\varepsilon_n}(\varphi_{\varepsilon_n}). \label{amontrer}
\end{eqnarray}

The full details of this fact can be found in  \cite[Lemma 3.1.]{bls}. We shall describe here only the  ideas of proof, which is achieved within two main steps. The first one consists in finding a bound $\Hh^1(K)\leq C$ when $\liminf  S_{\varepsilon_n}(\varphi_{\varepsilon_n})<+\infty $ (which is  obviously the case here).  

The main tool is the definition of the following geometric quantity: for each set $\Gamma\subset\R^2$, each unit vector $\nu\in\mathbb{S}^1$ and each $\lambda>0$ we set 
$$\Gamma_{\lambda,\nu}:=\{x\in\R^2\,:\,\mbox{ there exists  $t\in[-\lambda,\lambda]$ with } x-t\nu\in\Gamma\}$$
and we define
$$I_\lambda(\Gamma):=\frac{1}{2\pi\lambda} \int_{\mathbb{S}^1}\mathscr{(L}^2((\Gamma)_{\lambda,\nu}) d\nu.$$

The following geometrical estimate  \cite[Lemma 2.6.]{bls} is one of our key ingredients and  is of independent interest: whenever $\Gamma_\ve$ are compact connected sets converging to $\Gamma$ as $\ve\to 0$ in the Hausdorff distance, then 
\begin{equation}\label{geom est}
\exists \lambda,\varepsilon_0 >0 \; ; \quad I_\lambda(\Gamma_\varepsilon)\geq C \Hh^{1}(\Gamma_0), \quad \forall \varepsilon \leq \varepsilon_0, 
\end{equation}
where the constant $C$ is universal.

Now, fix $\delta_0,\tau_0>0$, and let $\{z_1, z_2, \dots, z_N\}\subseteq K$ be a $\tau_0$-network in $K$, i.e. $K\subseteq  \bigcup_{1\leq i\leq N} B(z_i,\tau_0)$.
Due to the convergence $d_{\varphi_\varepsilon}(z_i,x_\varepsilon)\to d(z_i)=0$, for small $\ve$ we can build a set $\Gamma_\ve=  \bigcup_{1\leq i \leq N} \Gamma_i^\varepsilon$ where each $\Gamma_i^\varepsilon$ is a $C^1$ curve connecting $z_i$ to  $x_1$ and satisfying $\int_{\Gamma_i^\varepsilon} \varphi_\varepsilon(s)d\Hh^1(s) < \delta_0$.

We use the usual estimate  $\frac{1}{4\varepsilon}(1-\varphi_\varepsilon)^2 + \varepsilon \|\nabla \varphi_\varepsilon\|^2 \geq \|\nabla (P(\varphi_\varepsilon))\|$ where $P(t)=t-t^2/2$ is a primitive of $(1-t)$, and compute the total variation of $P(\varphi_\varepsilon)$ in the direction $\nu$ on a set $(\Gamma_\ve)_{\lambda,\nu}$. Using that $P(\varphi_\varepsilon)$ is almost $0$ on $\Gamma_\ve$ (by definition of $\Gamma_\ve$ and using $P(t)\leq t$) and that, on the contrary, $P( \varphi_\varepsilon)\to P(1)=1/2$ a.e., we get an estimate on $I_\lambda(\Gamma_\ve)$. Thanks to \eqref{geom est} this turns into an estimate on the $\Hh^1$ measure of the Hausdorff limit of $\Gamma_\ve$. By taking then the limit $\delta_0\to 0$, and finally $\tau_0\to 0$ one gets an estimate on $\Hh^1(K)$ and concludes the first step.

The second step is a refinement of the first: once we have established the rectifiability of $K$, we can use the existence of tangent line $\Hh^1$-a.e. on $K$. Using a similar argument as the one above but adapted locally around each point of $K$ (i.e. choosing the direction $\nu$ orthogonal to the tangent to $K$ instead of taking an average over all directions) we are able to prove the better estimate \eqref{amontrer} and this finishes the proof. \qed






\begin{thebibliography}{00}


\bibitem{alr}
L.~Ambrosio, A.~Lemenant, and G.~Royer-Carfagni.
\newblock A variational model for plastic slip and its regularization via
  gamma-convergence.
\newblock {\em J. Elast., to appear}.

\bibitem{at}
L.~Ambrosio and V.~M. Tortorelli.
\newblock On the approximation of free discontinuity problems.
\newblock {\em Boll. Un. Mat. Ital. B (7)}, 6(1):105--123, 1992.

\bibitem{bcps}
F.~Benmansour, G.~Carlier, G.~Peyr{\'e}, and F.~Santambrogio.
\newblock Derivatives with respect to metrics and applications: subgradient
  marching algorithm.
\newblock {\em Numer. Math.}, 116(3):357--381, 2010.

\bibitem{bls}
M.~Bonnivard, A.~Lemenant, and F.~Santambrogio.
\newblock Approximation of length minimization problems among compact connected
  sets.
\newblock {\em in preparation}.

\bibitem{CarJimSan}
G.~Carlier, C.~Jimenez, and F.~Santambrogio.
\newblock Optimal transportation with traffic congestion and wardrop
  equilibria.
\newblock {\em SIAM J. Control Optim.}, 47:1330--1350, 2008.

\bibitem{gilbpoll}
E.~N. Gilbert and H.~O. Pollak.
\newblock Steiner minimal trees.
\newblock {\em SIAM J. Appl. Math.}, 16:1--29, 1968.

\bibitem{Karp}
R.~Karp.
\newblock Reducibility among combinatorial problems.
\newblock In {\em Complexity of Computer Computations}, pages 85--103. Plenum
  Press, 1972.

\bibitem{dmi}
G.~Dal Maso and F.~Iurlano.
\newblock Fracture models as {$\Gamma-$}limits of damage models.
\newblock {\em Comm. Pure Appl. Anal.}, 12(4):1657--1686., 2013.

\bibitem{MM}
L.~Modica and S.~Mortola.
\newblock Il limite nella {$\Gamma $}-convergenza di una famiglia di funzionali
  ellittici.
\newblock {\em Boll. Un. Mat. Ital. A (5)}, 14(3):526--529, 1977.

\bibitem{OudetKelvin}
{\'E}~Oudet.
\newblock Approximation of partitions of least perimeter by
  {$\Gamma$}-convergence: around {K}elvin's conjecture.
\newblock {\em Exp. Math.}, 20(3):260--270, 2011.

\bibitem{OudetSantam}
\'E. Oudet and F.~Santambrogio.
\newblock A modica-mortola approximation for branched transport and
  applications.
\newblock {\em Arch. Rati. Mech. An.}, 201(1):115--142, 2011.

\bibitem{ps2009}
E.~Paolini and E.~Stepanov.
\newblock Existence and regularity results for the steiner problem.
\newblock {\em Calc. Var. Partial Diff. Equations.}, 46(3):837--860, 2013.

\bibitem{sant}
Filippo Santambrogio.
\newblock A {M}odica-{M}ortola approximation for branched transport.
\newblock {\em C. R. Math. Acad. Sci. Paris}, 348(15-16):941--945, 2010.

\bibitem{sethian-fastmarching}
J.A. Sethian.
\newblock {\em Level Set Methods and Fast Marching Methods}.
\newblock Cambridge Monographs on Applied and Computational Mathematics.
  Cambridge University Press, 1999.




\end{thebibliography}

\end{document}